\documentclass[12pt]{amsart}
\usepackage{amscd,amssymb}
\usepackage[graph,frame,poly,arc]{xy}  
\usepackage[plainpages,backref,urlcolor=blue]{hyperref}

\theoremstyle{plain}
\newtheorem{thm}[subsection]{Theorem}
\newtheorem{lem}[subsection]{Lemma}

\newtheorem{cor}[subsection]{Corollary}

\theoremstyle{definition}

\newtheorem{ex}[subsection]{Example}

\numberwithin{equation}{section}
\newcommand{\A}{{\mathcal A}}

\newcommand{\C}{\mathbb{C}}
\newcommand{\PP}{\mathbb{P}}

\DeclareMathOperator{\im}{im}

\begin{document}

\title [On the global Tjurina numbers for line arrangements]
{On the minimal value of global Tjurina numbers for line arrangements}

\author[Alexandru Dimca]{Alexandru Dimca$^{1}$}
\address{Universit\'e C\^ ote d'Azur, CNRS, LJAD, France and Simion Stoilow Institute of Mathematics,
P.O. Box 1-764, RO-014700 Bucharest, Romania}
\email{dimca@unice.fr}

\thanks{$^1$ This work has been partially supported by the French government, through the $\rm UCA^{\rm JEDI}$ Investments in the Future project managed by the National Research Agency (ANR) with the reference number ANR-15-IDEX-01 and by the Romanian Ministry of Research and Innovation, CNCS - UEFISCDI, grant PN-III-P4-ID-PCE-2016-0030, within PNCDI III}



\begin{abstract} 
We show that a general lower bound for the global Tjurina number of a reduced complex projective plane curve, given by A. A. du Plessis and C.T.C. Wall, can be improved when the curve is a line arrangement.

\end{abstract}
 
\maketitle


\section{Introduction} 

Let $S=\C[x,y,z]$ be the graded polynomial ring in three variables $x,y,z$ and let $C:f=0$ be a reduced curve of degree $d$ in the complex projective plane $\PP^2$. The minimal degree of a Jacobian relation for $f$ is the integer $mdr(f)$
defined to be the smallest integer $m\geq 0$ such that there is a nontrivial relation
\begin{equation}
\label{rel_m}
 af_x+bf_y+cf_z=0
\end{equation}
among the partial derivatives $f_x, f_y$ and $f_z$ of $f$ with coefficients $a,b,c$ in $S_m$, the vector space of  homogeneous polynomials of degree $m$. We say that the plane curve has type $(d,r)$ if $d= \deg(f)$ and $r=mdr(f)$.

When $mdr(f)=0$, then $C$ is a pencil of lines, i.e. a union of lines passing through one point, a situation easy to analyze. We assume from now on that 
$$ mdr(f)\geq 1.$$
Denote by $\tau(C)$ the global Tjurina number of the curve $C$, which is the sum of the Tjurina numbers of the singular points of $C$. When $C$ is a line arrangement,  its global Tjurina number, coincides with its global Milnor number $\mu(C)$, and is given by
$$\tau(C)=\sum_p(n(p)-1)^2,$$
where the sum is over all multiple points $p$ of $C$ and $n(p)$ denotes the multiplicity of $C$ at $p$.
The following result is due to  A.A. du Plessis and C.T.C. Wall, see \cite[Theorem 3.2]{duPCTC}.
\begin{thm}
\label{thmCTC}

For a positive integer $r$, define two integers by 
$$\tau(d,r)_{min}=(d-1)(d-r-1) \text { and } \tau(d,r)_{max}=  (d-1)(d-r-1)+r^2.$$
For any reduced plane curve $C:f=0$ of type $(d,r)$, one has the inequalities
$$\tau(d,r)_{min} \leq \tau(C) \leq \tau(d,r)_{max}.$$
Moreover, suppose that  $r=mdr(f)>(d-1)/2$. Then one has the stronger inequality
$$\tau(C) \leq \tau(d,r)_{max}-{2r+2-d \choose 2}.$$
In particular, if $d$ is even and $r=d/2$, then
$ \tau(C) \leq \tau(d,r)_{max}-1.$
\end{thm}
It is easy to see that the lower bound is strict, i.e. there are curves $C$ of type $(d,r)$ such that
$\tau(C)=\tau(d,r)_{min}$ for any degree $d$ and $1 \leq r \leq d-1$, see \cite[Example 4.5]{DSt3syz} in general and \cite[Lemma 21]{E} for $(d,r)=(d,d-2)$.
The curves $C$ for which the maximal values for $\tau(C)$ are attained have special properties. Indeed, we have the following result, see \cite{Dmax}.

\begin{thm}
\label{thmDmax}
Let $C:f=0$ be a reduced plane curve of type $(d,r)$. Then the following hold.

\begin{itemize}
		\item[i)]  One has
$ \tau(C) =\tau(d,r)_{max}$
if and only if $C:f=0$ is a free curve with exponents $(r,d-1-r)$, and then $r<d/2$.
				\item[ii)]  One has
$ \tau(C) =\tau(d,r)_{max}-1$
if and only if $C:f=0$ is a nearly free curve with exponents $(r,d-r)$, and then $r \leq d/2$.
		
	\end{itemize}

\end{thm}
The free and nearly free curves, as well as related objects, are  actively investigated since some time, see \cite{ Abe18, AD, B+,Dcurves, Dmax,Drcc, DStExpo,DStFD,DStRIMS,DStMos,E,MaVa,ST,Yo}. Note that there are free  (resp. nearly free line arrangements)  in $\PP^2$ for any fixed type $(d,r)$ with $r<d/2$
(resp. with $r \leq d/2$), i.e. the maximal values for $\tau(C)$ can be attained with $C$ a line arrangement, at least when $r \leq d/2$.

In this note we consider the minimal values of $\tau(C)$ when $C:f=0$ is a {\it line arrangement}, once we fix its type $(d,r)$. Let $m(C)$ be the maximal multiplicity of a point in $C$,
and $n(C)$ the maximal multiplicity of a point in $C \setminus \{p\}$, where $p$ is any point in $C$ of multiplicity $m(C)$. Note that 
$$ 1 \leq n(C) \leq m(C) \leq d.$$
Moreover $m(C)=d$ if and only if $mdr(f)=0$, and $m(C)=d-1$ if and only if $mdr(f)=1$,
see \cite[Proposition 4.7]{DIM}. In addition, the case
$2=n(C)\leq m(C) \leq d-2$,  corresponds to the intersection lattice $L(C)$ being the lattice $L(d,m(C))$ discussed in \cite[Proposition 4.7]{DIM}, i.e. the intersection lattice of an arrangement obtained from $m(C)$ concurrent lines by adding $d-m(C)$ lines in general position.
The main result of this note is the following.

\begin{thm}
\label{thmA}
Let $C:f=0$ be an arrangement of $d\geq 4$ lines in $\PP^2$ which is not free. If we set $r=mdr(f)\geq 2 $, then  the following hold.

\begin{itemize}
		\item[i)]  With the above notation, one has
		$$\tau(C) \geq \tau'(d,r)_{min}:= \tau(d,r)_{min}+{r \choose 2}+{n(C) \choose 2}+1.$$
\item[ii)] If $r \ne d-m(C)$, then the possibly stronger inequality 	
$$\tau(C) \geq \tau''(d,r)_{min}:=\tau(d,r)_{min}+{r \choose 2} +{m(C) \choose 2}+1$$ holds.	

\end{itemize}

\end{thm}
The line arrangements such that $r=mdr(f) \in \{0,1,2\}$ are classified, see for instance \cite[Theorem 4.11]{DIM} for the case $r=2$, and for the remaining line arrangements we have the following.

\begin{cor}
\label{corA} Let $C:f=0$ be an arrangement of $d$ lines in $\PP^2$ which is not free and such that $r=mdr(f)\geq 3 $ and $n(C) \geq 3$. Then 
$$\tau(C) \geq \tau^N(d,r)_{min}:=\tau(d,r)_{min}+{r \choose 2} +4 \geq \tau(d,r)_{min}+7.$$
\end{cor}

\begin{ex}
\label{exA}
Consider the line arrangements $C:f=0$ of $d=7$ lines such that $r=mdr(f)=3$.
Then $\tau(7,3)_{min}=6\cdot 3=18$ and 
$$\tau'(7,3)_{min}=18+3+{n(C) \choose 2}+1=22+{n(C) \choose 2}.$$
The value $\tau=25$ is obtained for the following three line arrangements.
Two of them, say 
$$C_1:f_1=xyz(x+y)(x+3y)(x+2y+z)(4x+8y+z)=0$$ and 
$$C_2:f_2=xyz(x+2y+z)(y+z)(x-2y)(x-y)=0,$$
satisfy $m(C_i)=4$ and $n(C_i)=3$, since each of them has a quadruple point, a triple point and 12 double points. The two intersection lattices $L(C_1)$ and $L(C_2)$ are distinct, since only in $C_2$ there is a line in the arrangement containing both  points of multiplicity $>2$. Our lower  bound for both $C_1$ and $C_2$ is $\tau^N(7,3)_{min}=22+3=25$, i.e. a sharp lower bound.
The third line arrangement, say 
$$C_3:f_3=xyz(2x-3y+z)(x-y)(x+z)(y+z)=0,$$
 satisfy $m(C_3)=n(C_3)=3$, since it has 4 triple points and 9 double points. Hence we can use the bounds
$\tau'(d,r)_{min}=\tau''(d,r)_{min}=25$, since $3=mdr(C_3) \ne d-m(C_3)=4$. Therefore our various bounds are sharp in all these three situations.
\medskip

\end{ex}

\begin{ex}
\label{exA2}
Consider the line arrangements $C:f=0$ of $d$ lines, with intersection lattice $L(C)$ of type $\tilde L(m_1,m_2)$ as in \cite[Proposition 4.9]{DIM}, where $2 \leq m_1 \leq m_2$.
This line arrangement is the union of two pencils of lines, one consisting of $m_1$ lines, the other of $m_2$ lines, in general position to each other.
Then we have $d=m_1+m_2$, $r=mdr(f)=m_1$, $m(C)=m_2$, $n(C)=m_1$, 
$$\tau(C)=(d-1)^2-m_1m_2+1.$$
One has $\tau(d,r)_{min}=(d-1)^2-m_1m_2-(m_1^2-m_1)$,
while the lower bound given by Theorem \ref{thmA} i) is in this case
$$\tau'(d,r)_{min}=(d-1)^2-m_1m_2+1,$$
again a sharp bound.
\end{ex}

 Let $I_f$ denote the saturation of the ideal $J_f$ with respect to the maximal ideal ${\bf m}=(x,y,z)$ in $S$ and consider the local cohomology group 
 $$N(f)=I_f/J_f=H^0_{\bf m}(M(f)).$$
 The graded $S$-module  $N(f)$ satisfies a Lefschetz type property with respect to multiplication by generic linear forms, see  \cite{DPop}. This implies in particular the inequalities
$$0 \leq n(f)_0 \leq n(f)_1 \leq ...\leq n(f)_{[T/2]} \geq n(f)_{[T/2]+1} \geq ...\geq n(f)_T \geq 0,$$
where $T=3d-6$ and $n(f)_k=\dim N(f)_k$ for any integer $k$. We set
$$\nu(C)=\max _j \{n(f)_j\},$$
and call $\nu(C)$ the {\it freeness defect}  of the curve $C$. It is known that a curve $C$ is free (resp. nearly free) if and only if $\nu(C)=0$ (resp. $\nu(C)=1$), see \cite{Dmax}.
The relation between the invariants $\nu(C)$ and $mdr(f)$ is given by the following result, see \cite[Theorem 1.2]{Drcc}.
\begin{thm}
\label{thmN}
Let $C:f=0$ be a reduced plane curve of degree $d$ and let $r=mdr(f)$.
Then the following hold.
\begin{enumerate}
\item If $r < d/2$, then
$$\nu(C)=(d-1)^2-r(d-1-r)-\tau(C)=\tau(d,r)_{max}-\tau(C).$$

\item If $r \geq (d-2)/2$, then 
$$\nu(C)= \lceil   \  \frac{3}{4}(d-1)^2 \    \rceil  -\tau (C).$$

\end{enumerate}
\end{thm}
Here, for any real number $u$,  $\lceil     u    \rceil $ denotes the round up of $u$, namely the smallest integer $U$ such that $U \geq u$. Theorem \ref{thmA} yields the following.

\begin{cor}
\label{corB} 
Let $C:f=0$ be an arrangement of $d$ lines in $\PP^2$. If  we have $2 \leq r=mdr(f) <d/2$, then 
$$\nu(C) \leq \frac{r(r+1)}{2}-{n(C) \choose 2}-1 \leq \frac{r(r+1)}{2}-2.$$
\end{cor}

It is known that a line arrangement with $r=2$ is either free or nearly free, see \cite[Theorem 4.11]{DIM}. Hence in this case $\nu(C) \in \{0,1\}$, and the upper bound given by Corollary \ref{corB} is 1, hence this bound is sharp in this case. Another upper bound for the freeness defect $\nu(C)$ of a line arrangement is discussed in \cite{DStbound}.

\section{Proof of Theorem \ref{thmA} and of Corollary \ref{corB}} 
Consider the graded $S-$submodule $AR(f) \subset S^{3}$ of {\it all relations} involving the derivatives of $f$, namely
$$\rho=(a,b,c) \in AR(f)_q$$
if and only if  $af_x+bf_y+cf_z=0$ and $a,b,c$ are in $S_q$, the space of homogeneous polynomials of degree $q$. 
We have the following.
\begin{lem}
\label{L0} Let $C:f=0$ be a reduced plane curve of degree $d$, such that any irreducible component of $C$ is rational. Then
$$\tau(C)=\dim AR(f)_{d-2}+\frac{(d-1)(d-2)}{2}.$$
\end{lem}
\proof
Use \cite[Formula (3.4)]{DSe}
and the vanishing of $N(f)_{2d-3}$, to be defined below in the last section,  given in \cite[Theorem 2.7]{DStbound}.

\endproof

Let $\rho_1 \in AR(f)_r$ be a nontrivial relation of minimal degree. The inclusion
$S\cdot \rho_1 \subset AR(f)$ implies that 
$$\dim AR(f)_{d-2}=  \dim S_{d-r-2}+ \dim \overline {AR(f)}_{d-2},$$
where $\overline {AR(f)}=AR(f)/(S\rho_1)$. This implies
$$\tau(C)=\dim AR(f)_{d-2}+\frac{(d-1)(d-2)}{2}=(d-1)^2-r(d-\frac{r+1}{2})+\dim \overline {AR(f)}_{d-2}.$$
The exact sequence (3.3) in \cite{Dmax} shows that the graded $S$-module  $\overline {AR(f)}$ is torsion free. Let $p \in C$ (resp. $q \in C$) be two distinct points such that $p$ (resp. $q$) has multiplicity $m=m(C)$ (resp. $n=n(C)$).  Proceeding as in \cite[Section (2.2)]{Dcurves}, we get two syzygies, namely $\rho_p \in AR(f)_{d-m}$ and $\rho_q \in AR(f)_{d-n}$. We have the following.

\begin{lem}
\label{L1} Assume that $C$ is a line arrangement of $d$ lines.
Then  $\rho_q \notin S_{m-n} \cdot \rho_p$.
\end{lem}

\proof
Following \cite[Section (2.2)]{Dcurves}, we recall the construction of the two syzygies $\rho_p $ and $\rho_q $. First choose the coordinates on $\PP^2$ such that $p=(1:0:0)$ and $q=(0:1:0)$. To construct $\rho_p$, we write
$f=g(y,z)h(x,y,z)$, where $g=0$ is the equation of the $m=m(C)$ lines passing through $p$, and $h=0$ is the equation of the remaining $d-m$ lines in $C$. With this notation, one has
\begin{equation}
\label{eqA}
\rho_p=(xh_x-d\cdot h,yh_x,zh_x),
\end{equation}
where $h_x$ denotes the partial derivative of $h$ with respect to $x$. Moreover, one has
$f_x=gh_x$. Similarly, one can write $f=\tilde g(x,z)\tilde h(x,y,z)$, where $\tilde g=0$ is the equation of the $m=n(C)$ lines passing through $q$, and $\tilde h=0$ is the equation of the remaining $d-m$ lines in $C$. With this notation, one has
\begin{equation}
\label{eqB}
\rho_q=(x\tilde h_y,y\tilde h_y-d\cdot \tilde h,z\tilde h_y),
\end{equation}
where $\tilde h_y$ denotes the partial derivative of $\tilde h$ with respect to $y$. Moreover, one has
$f_y=\tilde g \tilde h_y$.
Assume there is a polynomial $A \in S_{m-n}$ such that $\rho_q =A \rho_p$. Looking at the third coordinate, we get $\tilde h_y=Ah_x$.
Then the equality on the first coordinate implies $Ah=0$, which is a contradiction.
\endproof

Now we complete the proof of Theorem \ref{thmA}, and in view of Example \ref{exA2} we may assume that $C$ is not of type $\tilde L(m_1,m_2)$.
We apply \cite[Theorem 1.2]{Dcurves} and discuss the following cases, in order to prove the claim (i).

\medskip

\noindent {\bf Case 1}: $r=mdr(f)=d-m$. Then we can take $\rho_1=\rho_p$ and the class of $\rho_q$ in $\overline {AR(f)}$ is nonzero, because of Lemma \ref{L1} and of the fact that the syzygies constructed in \cite[Section (2.2)]{Dcurves} are primitive, i.e. they are not multiple of strictly lower degree syzygies. It follows that the vector space
$S_{n-2}\rho_q$ is naturally embedded in $\overline {AR(f)}_{d-2}$.

We recall now the construction of the Bourbaki ideal $B(C,\rho_1)$ associated to a degree $d$ reduced curve $C:f=0$ and to a minimal degree non-zero syzygy $\rho_1 \in AR(f)$, see \cite{DStJump}.
For any choice of the syzygy $\rho_1=(a_1,b_1,c_1)\in AR(f) $ with minimal degree $r=d_1$, we have a morphism of graded $S$-modules
\begin{equation} \label{B1}
S(-r)  \xrightarrow{u} AR(f), \  u(h)= h \cdot \rho_1.
\end{equation}
For any homogeneous syzygy $\rho=(a,b,c) \in AR(f)$, consider the determinant $\Delta(\rho)=\det M(\rho)$ of the $3 \times 3$ matrix $M(\rho)$ which has as first row $x,y,z$, as second row $a_1,b_1,c_1$ and as third row $a,b,c$. Then it turns out that $\Delta(\rho)$ is divisible by $f$, see \cite{Dmax}, and we define thus a new morphism of graded $S$-modules
\begin{equation} \label{B2}
 AR(f)  \xrightarrow{v}  S(r-d+1)   , \  v(\rho)= \Delta(\rho)/f,
\end{equation}
and a homogeneous ideal $B(C,\rho_1) \subset S$ such that $\im v=B(C,\rho_1)(r-d+1)$.
It is known that the ideal $B(C,\rho_1)$, when $C$ is not a free curve, defines a $0$-dimensional subscheme $Z(C,\rho_1)$ in $\PP^2$, which is locally a complete intersection, see \cite[Theorem 5.1]{DStJump}.
This construction yields the following exact sequence.
$$0 \to S(m-d)  \xrightarrow{u} AR(f) \xrightarrow{v}  B(C,\rho_1)(-m+1) \to  0.$$
Hence, with this notation,  $h_q=v(\rho_q) \in B_{d-n-m+1}$ and
$$S_{n-2}h_q  \subset \overline {AR(f)}_{d-2}  = B_{d-m-1}.$$
On the other hand we know, see \cite[Corollary 3.5]{SchCMH} or \cite[Corollary 3.6]{DIM}
, that  the Castelnuovo-Mumford regularity $reg(AR(f))$ of the graded $S$-module $AR(f)$ satisfies $reg(AR(f))\leq d-2$ as soon as $d \geq 4$. In particular, all the generators in a minimal set of generators for $AR(f)$ have degrees $\leq d-2$.
It follows that our ideal $B$ is generated in degrees $\leq d-m-1$. Moreover we know that the subscheme in $\PP^2$ defined by the ideal $B$ is 0-dimensional, \cite[Theorem 4.1]{DStJump}, and hence $B$ cannot be generated by a single polynomial $h_q$. It follows that there is at least one polynomial $\hat h \in B_{d-m-1} \setminus S_{n-2}h_q$ and hence
$$\dim AR(f)_{d-2} \geq \dim S_{n-2} +1.$$
This completes the proof of the claim $(i)$ in this case.
\medskip

\noindent {\bf Case 2}: $r=mdr(f)=m-1$ and $C$ is free, case which is discarded by our hypotheses.

\medskip

\noindent {\bf Case 3}: $m\leq r=mdr(f) \leq d-m-1$. Then let $\rho_1 \in AR(f)_r$ by a minimal degree non zero syzygy. Exactly as in the Case 1 above,
it follows that the vector space
$S_{m-2}\rho_p$ is naturally embedded in $\overline {AR(f)}_{d-2}=B_{r-m-1}$, and this proves the claims (i) and (ii) in this case as well, exactly as above.

\subsection{Proof of Corollary \ref{corB}} 

 One has
$$\nu(C) = \tau(d,r)_{max}-\tau(C) \leq \tau(d,r)_{max}-\tau'(d,r)_{min}=$$
$$=\frac{r(r+1)}{2}-{n(C) \choose 2}-1\leq \frac{r(r+1)}{2}-2.$$

\section{An application to Terao's conjecture}

H. Terao has conjectured that if $\A$ and $\A'$ are hyperplane arrangements in $\PP^n$  with isomorphic intersection lattices $L(\A) =L(\A')$, and if $\A$ is free, then $\A'$ is also free, see  for details
\cite{DHA, OT, Yo}.
Using Theorem \ref{thmA}, we get the following partial positive answer in the case of line arrangements, which is an improvement of \cite[Corollary 2.3]{Dmax}.

\begin{cor}
\label{cor2} Let $C:f=0$ and $C':f'=0$ be two arrangements of $d$ lines  in $\PP^2$  with isomorphic intersection lattices $L(C) =L(C')$.  If $C$ is free, then there is  a unique integer $r\geq 0$ such that $r<d/2$ and
$$\tau(\A) =(d-1)(d-r-1)+r^2.$$
If this integer $r$ satisfies 
$$r<\frac{-3+\sqrt {8d+41}}{2},$$
 then the line arrangement $C'$ is also free.

\end{cor}

\proof We can assume that $n(C) \geq 3$, since in the other cases Terao's conjecture is obvious. Indeed, for $n(C) \leq 2$, the intersection lattice $L(C)$ has the type $L(d,m(C))$, and the freeness of $C$ depends only on $m(C)$, see \cite[Proposition 4.7]{DIM}.
When $n(C) \geq 3$, a direct  computation shows that for $s<\frac{-3+\sqrt {8d+41}}{2}$, the intervals
$[\tau^N(d,s)_{min},\tau(d,s)_{max}]$ and $[\tau^N(d,s-1)_{min},\tau(d,s-1)_{max}]$ are disjoint, since
$$\tau(d,s)_{max} <\tau^N(d,s-1)_{min}.$$
It follows that each value $\tau(d,s)_{max}$ uniquely determines the corresponding $s$ when $s<\frac{-3+\sqrt {8d+41}}{2}$.
Moreover, $\tau(C)=\tau(C')$, since this number depends only on the intersection lattice  $L(C) =L(C')$. It follows that $mdr(f')=mdr(f)=r$ and hence $C'$ is free by Theorem  \ref{thmDmax} $i)$.

\endproof

\begin{ex}
\label{exE} When $d=100$, then we have
$$\frac{-3+\sqrt {8d+41}}{2}=\frac{-3+29}{2}=13.$$
Hence Terao's conjecture holds for any free arrangement $C:f=0$ of 100 lines with $r\leq 13$. The old result  \cite[Corollary 2.3]{Dmax} when $d=100$ gives the bound $r \leq \sqrt {d-2}$, hence $ r \leq 9$.

\end{ex}



\begin{thebibliography}{00}


\bibitem{Abe18}  T. Abe, Plus-one generated and next to free arrangements of hyperplanes, arXiv:1808.04697.

\bibitem{AD}  T. Abe, A. Dimca, On the splitting types of bundles of logarithmic vector fields along plane curves, Internat. J. Math. 29 (2018), no. 8, 1850055, 20 pp.

\bibitem{B+} E. Artal Bartolo, L. Gorrochategui, I. Luengo, A. Melle-Hern\' andez,
On some conjectures about free and nearly free divisors, in: {\it Singularities and Computer Algebra, Festschrift for Gert-Martin Greuel on the Occasion of his 70th Birthday}, pp. 1--19, Springer (2017)









\bibitem{Dcurves}  A. Dimca, Curve arrangements, pencils, and Jacobian syzygies,  Michigan Math. J. 66 (2017), 347--365.

\bibitem{DHA}  A. Dimca,   {\em Hyperplane Arrangements: An Introduction}, Universitext, Springer, 2017.


\bibitem{Dmax}  A. Dimca, Freeness versus maximal global Tjurina number for plane curves, 
Math. Proc. Cambridge Phil. Soc.  163 (2017), 161--172.





\bibitem{Drcc}  A. Dimca, On rational cuspidal plane curves, and the local cohomology of Jacobian rings, 
arXiv:1707.05258. to appear in Commentarii Mathematici Helvetici.

\bibitem{DIM} A. Dimca, D. Ibadula, A. M\u acinic, Numerical invariants and moduli spaces for line arrangements, arXiv:1609.06551.

\bibitem{DPop} A. Dimca, D. Popescu, 
Hilbert series and Lefschetz properties of dimension one almost complete intersections, Comm. Algebra 44 (2016), 4467--4482.


\bibitem{DSe} A. Dimca, E. Sernesi,  Syzygies and logarithmic vector fields along plane curves, Journal de l'\'Ecole polytechnique-Math\'ematiques 1(2014), 247-267.


\bibitem{DStExpo} A. Dimca, G. Sticlaru, On the exponents of free and nearly free projective plane curves, Rev. Mat. Complut. 30(2017), 259--268.

\bibitem{DStFD} A. Dimca, G. Sticlaru, Free divisors and rational cuspidal plane curves, Math. Res. Lett. 24(2017), 1023--1042.



\bibitem{DStRIMS} A. Dimca, G. Sticlaru, Free and nearly free curves vs. rational cuspidal plane curves, Publ. RIMS Kyoto Univ. 54 (2018), 163--179.

\bibitem{DStMos} A. Dimca, G. Sticlaru, On the freeness of rational cuspidal plane curves, Moscow Math. J. 
18(2018), 659--666.


\bibitem{DStbound} A. Dimca, G. Sticlaru, Line and rational curve arrangements, and Walther's inequality, arXiv:1803.05386.

\bibitem{DStJump} A. Dimca, G. Sticlaru, On the jumping lines of bundles of logarithmic vector fields along plane curves, arXiv: 1804.06349.

\bibitem{DSt3syz} A. Dimca, G. Sticlaru, Plane curves with three syzygies, plus-one generated curves, and nearly cuspidal curves, arXiv: 1810.11766.



\bibitem{duPCTC} A.A. du Pleseis,  C.T.C. Wall, Application of the theory of the
discriminant to highly singular plane curves, Math. Proc. Camb.
Phil. Soc.,  126(1999), 259-266. 


\bibitem{E} Ph. Ellia,  
Quasi complete intersections and global Tjurina number of plane curves, arXiv:1901.00809.

 \bibitem{MaVa} S. Marchesi, J. Vall\` es, Nearly free curves and arrangements: a vector bundle point of view, arXiv:1712.04867.


\bibitem{OT}
P. Orlik and H. Terao,
{\em Arrangements of Hyperplanes,}
Springer-Verlag, Berlin Heidelberg New York, 1992.





\bibitem{SchCMH} H. K. Schenck, Elementary modifications and line configurations in $\PP^2$,  Comment. Math. Helv. 78 (2003), 447--462.


\bibitem{ST} A. Simis, S.O. Toh\u aneanu, Homology of homogeneous divisors, Israel J. Math. 200 (2014), 449-487.




\bibitem{Yo} M. Yoshinaga, Freeness of hyperplane arrangements and related topics, Annales de la Facult\'e des Sciences de Toulouse, vol. 23 no. 2 (2014), 483-512.

\end{thebibliography}
\end{document}